\documentclass[11pt]{article}

\usepackage{amsmath,amsfonts,amssymb,amsthm,graphicx}

\begin{document}

\newtheorem{proposition}{Proposition}[section]
\newtheorem{theorem}{Theorem}[section]
\newtheorem{remark}{Remark}[section]
\newtheorem{example}{Example}[section]
\newtheorem{corollary}[theorem]{Corollary}
\newtheorem{lemma}{Lemma}[section]
\newtheorem{definition}{Definition}[section]

\baselineskip3.4ex

\def\arraystretch{1.75}

\allowdisplaybreaks

\title{On convergence of importance sampling and other properly weighted samples to the target distribution}
\author{Sonia Malefaki and George Iliopoulos\footnote{Corresponding author; Address: Department of Statistics and Insurance Science, University of Piraeus, 80 Karaoli \& Dimitriou str., 18534 Piraeus, Greece, e--mail: {\tt geh@unipi.gr} }}

\maketitle

\begin{abstract} We consider importance sampling as well as other properly weighted samples with respect to a target distribution $\pi$ from a different point of view. By considering the associated weights as sojourn times until the next jump, we define appropriate jump processes. When the original sample sequence forms an ergodic Markov chain, the associated jump process is an ergodic semi--Markov process with stationary distribution $\pi$. Hence, the type of convergence of properly weighted samples may be stronger than that of weighted means. In particular, when the samples are independent and the mean weight is bounded above, we describe a slight modification in order to achieve exact (weighted) samples from the target distribution. 
\end{abstract}




\section{Introduction}

Importance sampling (IS) (Marshall, 1956) is a well--known Monte Carlo method that is useful in any discipline where integral approximations are needed. For a measurable space $(\mathcal X,\mathcal B(\mathcal X))$ and  a probability distribution $\pi$ defined on it, the method attempts to estimate the integral
\[
\mathbf E_{\pi}(h) := \int_{\mathcal X} h(x)\pi({\rm d}x), 
\]
for $h\in\mathcal L^{1}(\pi)$, by drawing independent samples $x_{1},\ldots,x_{n}$ from a trial distribution $g$ with support at least that of $\pi$. Assuming that the distributions have densities $\pi(x)$, $g(x)$ with respect to a $\sigma$--finite measure $\mu$ [i.e., $\pi({\rm d}x)=\pi(x)\mu({\rm d}x)$, $g({\rm d}x)=g(x)\mu({\rm d}x)$], $\mathbf E_{\pi}(h)$ is approximated by
\begin{equation} \label{IS estimate n}
\bar{h}_{n}^{IS} := \frac{\sum_{i=1}^{n}w(x_{i})h(x_{i})}{n},
\end{equation}
or by
\begin{equation} \label{IS estimate sumw}
\hat{h}_{n}^{IS} := \frac{\sum_{i=1}^{n}w(x_{i})h(x_{i})}{\sum_{i=1}^{n}w(x_{i})},
\end{equation}
where $w(x_{i}):=\pi(x_{i})/g(x_{i})$. In standard terminology, $g$ is called ``the importance distribution'' and $w(x_{i})$'s ``the importance weights''. 

The validity of $\bar{h}_{n}$ and $\hat{h}_{n}$ as approximations of $\mathbf E_{\pi}(h)$ is justified by the strong law of large numbers which ensures that for all $h\in\mathcal L^{1}(\pi)$ it holds $n^{-1}\sum w(X_{i})h(X_{i})\rightarrow\mathbf E_{g}\{w(X)h(X)\}=\mathbf E_{\pi}(h)$ and $n^{-1}\sum w(X_{i})\rightarrow\mathbf E_{g}\{w(X)\}=1$ with probability one. Note that $\hat{h}_{n}$ can be used in more general settings than $\bar{h}_{n}$, e.g.\ when the importance weights $w(x)$ are known up to a multiplicative constant. In this case the second limit in general differs from one whilst the first becomes $\mathbf E_{g}\{w(X)\}\mathbf E_{\pi}(h)$. Also note that the assumption of independent samples can be relaxed since there are more general contexts under which the above IS estimators converge to $\mathbf E_{\pi}(h)$. For example, if the sequence of $X$'s forms a Harris ergodic Markov chain with stationary distribution $g$, the above limits still hold due to the Ergodic Theorem. 

Besides integral approximations, simulation methods aim to generate samples from a target distribution. In this respect, IS seems at first glance to fail obtaining samples from $\pi$ since all draws are made from the importance distribution $g$. The Sampling/Importance Resampling method (SIR, Rubin, 1987) is an attempt to circumvent this drawback by weighted resampling with replacement from the generated $g$--sample. The weight assigned to $x_{i}$ is its normalised importance weight $w(x_{i})/\sum_{i=1}^{n}w(x_{i})$. As $n\rightarrow\infty$, this approach produces a sample which is approximately $\pi$--distributed. Smith and Gelfand (1992) revisited the approach and proved the assertion by considering the resampling procedure as weighted bootstrapping. 

The aim of the present paper is to show that, under certain conditions, when the $g$--sample is properly weighted it converges in a sense to the target distribution $\pi$. Actually, this is true for a jump process associated with the weighted sample. From this point of view, importance weighting does not differ much from MCMC sampling schemes. In fact, some of them are special cases of the above mentioned jump processes (e.g.\ the Metropolis--Hastings algorithm, see Subsection 3.1). It turns out that in order to obtain approximate samples from $\pi$ resampling is not needed at all. However, our intention is to present these facts without to criticise SIR. 

The paper is organised as follows: In Section 2 we define the jump process associated with a weighted sample and prove that under certain conditions it converges weakly to the target distribution. Section 3 contains some examples of known sampling schemes which are special cases of this context. As a result of independent interest, we give an upper bound for the total variation distance to stationarity of Markov jump processes with embedded Doeblin chains when the mean sojourn time is bounded above. In Section 4 we discuss the case of stationary weighted sequences and show that in some cases it is possible to locate the time after which the associated jump process reaches equilibrium. Finally, an Appendix contains proofs of the stated propositions.


\section{Jump processes associated with properly weighted sequences}

The concept of a properly weighted sample has been introduced by Liu and Chen (1998) (see also Liu, 2001) as a generalisation of the standard IS method. 

\begin{definition} \label{properly weighted sample} {\rm [Liu and Chen (1998)]} A set of weighted random samples $(X_{i},\xi_{i})_{1\leqslant i\leqslant n}$ is called proper with respect to $\pi$ if for any square integrable function $h$, 
\[
\mathbf E\{\xi_{i}h(X_{i})\} = \kappa \mathbf E_{\pi}\{h(X_{i})\},\ \mbox{for}\ i=1,\ldots,n,
\]
for some positive constant $\kappa$. 
\end{definition}

As Liu (2001) points out, an equivalent definition is the following:

\noindent
{\bf Definition \ref{properly weighted sample}(a)} \label{properly weighted sample a}
{\it
A set of weighted random samples $(X_{i},\xi_{i})_{1\leqslant i\leqslant n}$ is called proper with respect to $\pi$ if 
\[
\mathbf E\{\xi_{i}|X_{i}=x\} = \kappa\pi(x)/g(x),\ \mbox{for}\ i=1,\ldots,n,
\]
for some positive constant $\kappa$, where $X_{i}\sim g$.
}

In the sequel, we will associate with any infinite weighted random sequence a jump process in the following sense.  

\begin{definition} \label{properly weighted sequence}
Consider a weighted sequence $(X_{n},\xi_{n})_{n\in\mathbf Z_{+}}:=((X_{0},\xi_{0}),(X_{1},\xi_{1}),\ldots)$, where the $\xi$'s are strictly positive weights. Define $S_{0}=0$, $S_{n}=\sum_{i=0}^{n-1}\xi_{i}$, $n\geqslant 1$, and let 
\[
N_{t}:=\sup\{n: S_{n}\leqslant t\},\ t\geqslant 0.
\]
Then, the stochastic process $Y=(Y_{t})_{t\geqslant 0}$ defined by $Y_{t}:=X_{N(t)}$, $t\geqslant 0$, will be called the jump process associated with the weighted sequence $(X_{n},\xi_{n})_{n\in\mathbf Z_{+}}$. 
\end{definition}

The definition ensures that the process $Y$ has right continuous sample paths which also have left hand limits. However, if the support of $\xi_{n}$'s is a subset of $\mathbf N=\{1,2,\ldots\}$, we will consider the process $Y$ only for $t\in\mathbf Z_{+}$, i.e.\ we set $Y=(Y_{0},Y_{1},Y_{2},\ldots)$. If this is the case, limits of quantities related to $Y_{t}$ should be suitably interpreted. 

\begin{proposition} \label{lim Y equals pi}
Assume that the sequence $X=(X_{n})_{n\in\mathbf Z_{+}}$ is a homogeneous Harris ergodic Markov chain with state space $(\mathcal X,\mathcal B(\mathcal X))$ having an invariant probability distribution $g$ and the distribution of $\xi_{n}$ depends solely on $X_{n}$ with $\mathbf E\{\xi_{n}|X_{n}=x\} = \kappa w(x) = \kappa\pi(x)/g(x)$ for some $\kappa>0$. Then, for the jump process $(Y_{t})_{t\geqslant 0}$ associated with the weighted sequence $(X_{n},\xi_{n})_{n\in\mathbf Z_{+}}$ it holds that
\[
\lim_{t\uparrow\infty}\mathbf P\{Y_{t}\in A\} = \pi(A),\ \forall\:A\in\mathcal B(\mathcal X).
\]
\end{proposition}
\begin{proof}
The result follows from the standard theory of semi--Markov processes (cf.\ Limnios and Opri\c{s}an, 2001). Under the above assumptions, $Y$ is an ergodic semi--Markov process with embedded Markov chain $X$ and respective sojourn times $(\xi_{n})_{n\in\mathbf Z_{+}}$. Thus, 
\[
\lim_{t\uparrow\infty}\mathbf P\{Y_{t}\in A\}
= 
\frac{\int_{A} \mathbf E\{\xi|x\} g(x)\mu({\rm d}x)}{\int_{\mathcal X} \mathbf E\{\xi|x\} g(x)\mu({\rm d}x)} 
=
\frac{\int_{A} \kappa w(x) g(x)\mu({\rm d}x)}{\int_{\mathcal X} \kappa w(x) g(x)\mu({\rm d}x)} 
=
\frac{\int_{A} \pi(x)\mu({\rm d}x)}{\int_{\mathcal X} \pi(x)\mu({\rm d}x)} 
= 
\pi(A)
\]
as is claimed.
\end{proof}

Setting deterministically $\xi_{n}\equiv w(X_{n})$, we have the following:

\begin{corollary} \label{corollary: IS}
If $(X_{n})_{n\in\mathbf Z_{+}}$ forms a Harris ergodic Markov chain with stationary distribution $g$, then the jump process associated with the weighted sequence $(X_{n},w(X_{n}))_{n\in\mathbf Z_{+}}$ has $\pi$ as limit distribution.
\end{corollary}

Any sequence of independent $g$--distributed random variables forms trivially an ergodic Markov chain with stationary distribution $g$. Thus, Corollary \ref{corollary: IS} covers also the original importance weighted sequence. 

\begin{figure}[t]
\begin{small}
\centerline{
\begin{tabular}{ccc}
$Y_{1}$ & $Y_{3}$ & $Y_{10}$ \\
\includegraphics[scale=0.45]{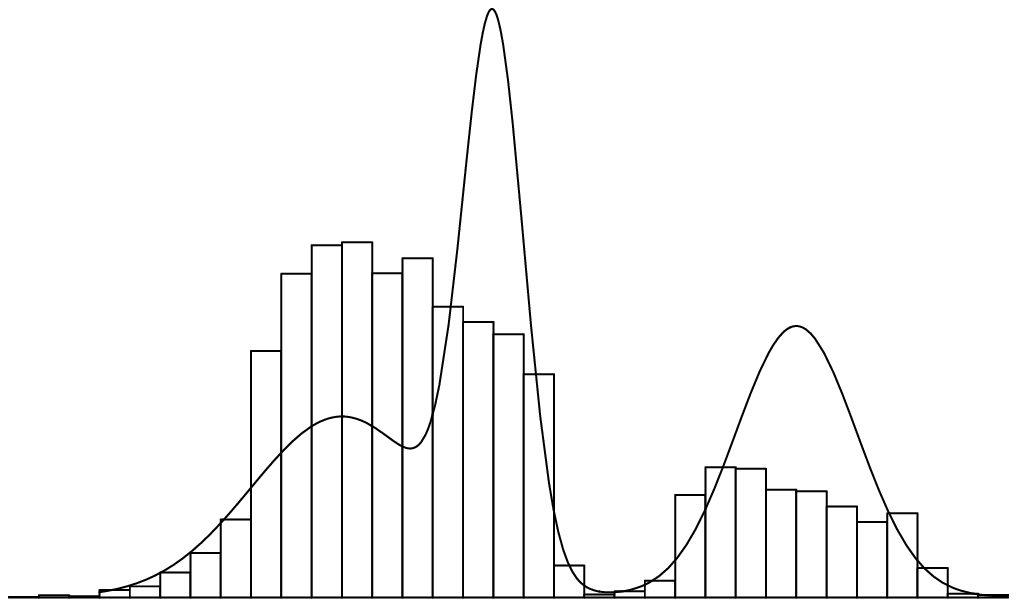} & \includegraphics[scale=0.45]{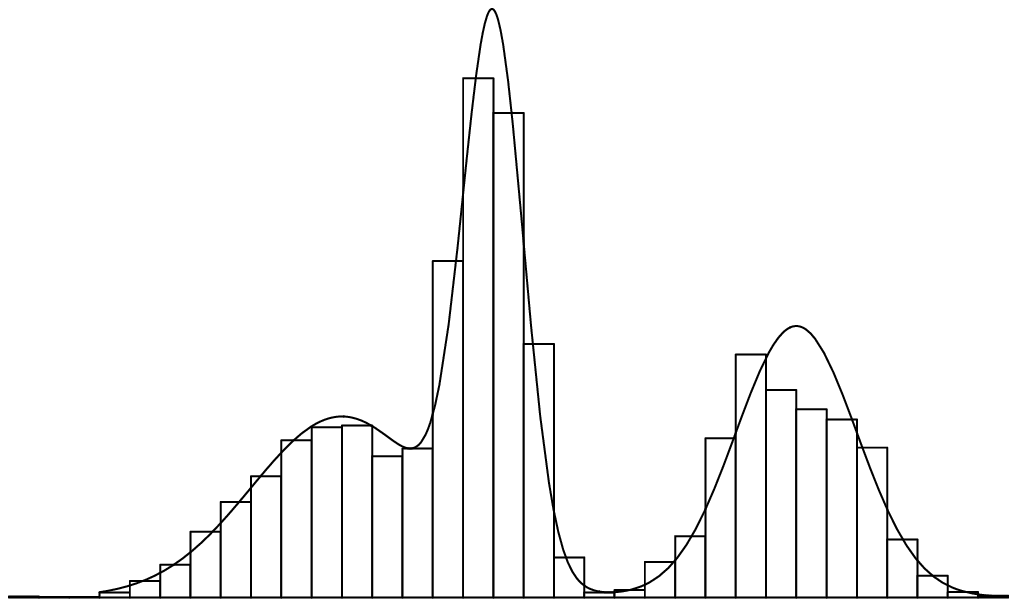} & \includegraphics[scale=0.45]{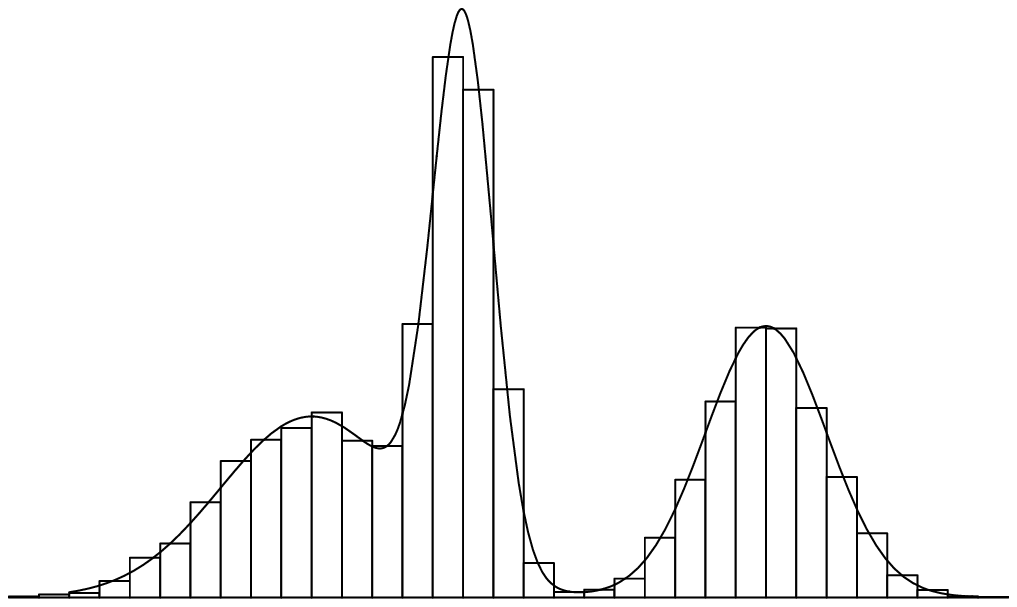}
\end{tabular}
}
\caption{Histograms of the values of $Y_{1}$, $Y_{3}$ and $Y_{10}$ of $m=10000$ independent IS runs with target distribution $\pi\sim \frac{1}{3}\:\mathcal N(0,9)+\frac{1}{3}\:\mathcal N(5,1)+\frac{1}{3}\:\mathcal N(15,4)$ and $g\sim\mathcal C(0,10)$.}
\label{figure:cauchy for normal mixture}
\end{small}

\bigskip
\hrule
\end{figure}

\begin{example} \label{IS normal mixture}
{\rm
Let the target distribution be the normal mixture $\pi\sim\frac{1}{3}\:\mathcal N(0,3^{2})+\frac{1}{3}\:\mathcal N(5,1)+\frac{1}{3}\:\mathcal N(15,2^{2})$ and $g\sim\mathcal C(0,10)$ (i.e.\ centered Cauchy with scale $10$). We have run the standard IS algorithm $m=10000$ times independently and for each run we have recorded the values $Y_{1}$, $Y_{3}$ and  $Y_{10}$ of the associated jump process. The corresponding histograms in Figure \ref{figure:cauchy for normal mixture} clearly illustrate the distributional convergence to $\pi$.
}

\end{example}

Under the assumptions of Proposition \ref{lim Y equals pi}, the weighted average 
\begin{equation} \label{xi estimate}
\hat{h}_{n} := \frac{\sum_{i=0}^{n-1}\xi_{i}h(X_{i})}{\sum_{i=0}^{n-1}\xi_{i}}
\end{equation}
converges almost surely to $\mathbf E_{\pi}(h)$ for any $h\in\mathcal L^{1}(\pi)$. Furthermore, if the Central Limit Theorem holds for $\hat{h}_{n}$, then its asymptotic variance is
\[
\sigma^{2}(h) = \sigma^{2}_{IS}(h) + \frac{1}{\kappa^{2}}\:\mathbf E_{g}\!\left[{\rm Var}\{\xi|X\}\{h(X)-\mathbf E_{\pi}(h)\}^{2}\right],  
\]
where $\sigma^{2}_{IS}(h)$ is the asymptotic variance of the IS estimator $\hat{h}_{n}^{IS}$ in (\ref{IS estimate sumw}) (which arises using $\xi_{n}\equiv \kappa w(X_{n})$ deterministically). This follows immediately from the identity
\[
\mathbf E_{g}\{\xi_{i}h(X_{i})\xi_{j}f(X_{i})\} = 
\left\{
\begin{array}{ll}
\kappa^{2}\mathbf E_{g}\{w(X_{i})^{2}h(X_{i})f(X_{j})\} + \mathbf E_{g}\{{\rm Var}(\xi_{i}|X_{i})h(X_{i})f(X_{i})\},&\hspace*{-1ex}i= j,
\\
\kappa^{2}\mathbf E_{g}\{w(X_{i})h(X_{i})w(X_{j})f(X_{j})\},&\hspace*{-1ex}i\neq j,
\end{array}
\right.
\]
and the formula for the asymptotic variance of a ratio estimator. As is intuitively rational, randomisation of the weights increases the variance. Hence, for estimation purposes the IS estimator is always more accurate. Otherwise, the less variable the weights are the more accuracy is achieved.

In light of Definition \ref{properly weighted sequence}, the estimator (\ref{xi estimate}) can be expressed as
\[
\hat{h}_{n} = \frac{1}{S_{n}}\int_{0}^{S_{n}}h(Y_{s})\nu({\rm d}s), 
\]
where $\nu$ is either the Lebesgue or the counting measure depending on whether the weights are continuous or discrete. Note that in the case of standard IS estimator (i.e.\ when $\xi_{n}\equiv \kappa w(X_{n})$), this gives another justification for the use of $\hat{h}_{n}^{IS}$ in (\ref {IS estimate sumw}) instead of $\bar{h}_{n}^{IS}$ in (\ref{IS estimate n}). 

Proposition \ref{lim Y equals pi} establishes a stronger result than that of the convergence of Ces\`aro averages
\[
\hat{h}_{t}:=t^{-1}\int_{0}^{t}h(Y_{s})\nu({\rm d}s)
\] 
to $\mathbf E_{\pi}(h)$. It states that there is distributional convergence of the generated sequence to the target distribution analogous to that of MCMC schemes. Thus, properly weighted samples and in particular the output of the standard IS method are in general serious competitors of them. 

The requirement that the distribution of $\xi_{n}$ depends only on $x_{n}$ seems rather restrictive. However, $x_{n}$ could be a block of specific size allowing $\xi_{n}$ to depend on more than one term of the sequence. (Note that the standard definition of a semi--Markov process allows the sojourn time of $X_{n}$ depending on both $X_{n}$ and $X_{n+1}$.)
This is illustrated via two examples in Subsection 3.4. 


\section{Examples}

In this section we discuss some known simulation schemes which are special cases of the jump process context. More specifically, we refer to cases where the conditional distribution of the weights is geometric (discrete case) or exponential (continuous case) and thus the associate jump process is a pure Markov jump process. Moreover, we consider two IS estimators used in diagnosing convergence of MCMC schemes under the current perspective. 

\subsection{The Metropolis--Hastings algorithm}

Consider an arbitrary Metropolis--Hastings (MH, Metropolis et al., 1953; Hastings, 1970) algorithm with target distribution $\pi$ and proposal $q(\cdot|\cdot)$, that is, at time $t+1$ given $Y_{t}=y$, draw $Z\sim q(z|y)$ and set $Y_{t+1}=z$ with probability 
\begin{equation} \label{MH acceptance probability}
a(y,z)=\min\left\{1,\frac{\pi(z)q(y|z)}{\pi(y)q(z|y)}\right\}
\end{equation}
or $Y_{t+1}=y$ with probability $1-a(y,z)$. Although it is well--known that the algorithm defines a reversible Markov chain with stationary distribution $\pi$, let us consider it from a different point of view. 

Let $X=(X_{n})_{n\in\mathbf Z_{+}}$ be a Markov chain with transition density
\[
g(x_{i}|x_{i-1})=\frac{a(x_{i-1},x_{i})q(x_{i}|x_{i-1})}{\int a(x_{i-1},z)q(z|x_{i-1})\mu({\rm d}z)} = 
\frac{\min\{\pi(x_{i-1})q(x_{i}|x_{i-1}),\pi(x_{i})q(x_{i-1}|x_{i})\}} {\int\min\{\pi(x_{i-1})q(z|x_{i-1}),\pi(z)q(x_{i-1}|z)\}\mu({\rm d}z)}.
\]
(Notice that this is exactly the density of the accepted states of the above MH algorithm.) It can be easily verified that $g(x_{i}|x_{i-1})$ satisfies the detailed balance condition
\[
g(x_{i-1})g(x_{i}|x_{i-1}) = g(x_{i})g(x_{i-1}|x_{i}),
\]
where $g(x)\propto \int \min\{ \pi(x)q(z|x),\pi(z)q(x|z)\}\mu({\rm d}z)$. This function, when normalised, results to a probability density function since 
\[
\textstyle\int g(x)\mu({\rm d}x) 
\leqslant \int \int \pi(x)q(z|x)\mu({\rm d}z)\mu({\rm d}x) = 1,  
\]
hence it is the stationary distribution of the Markov chain $X$.  Weight now $x_{i}$ by $\xi_{i}$ drawn from the geometric distribution with probability mass function
\[
p(\xi|x_{i}) = \left\{\textstyle\int a(x_{i},z)q(z|x_{i})\mu({\rm d}z)\right\}
\left\{1-\textstyle\int a(x_{i},z)q(z|x_{i})\mu({\rm d}z)\right\}^{\xi-1},\ \xi=1,2,\ldots.
\]
Since
\begin{equation} \label{MH E xi}
\mathbf E\{\xi|x_{i}\} = \left\{\textstyle\int a(x_{i},z)q(z|x_{i})\mu({\rm d}z)\right\}^{-1} \propto \pi(x_{i})/g(x_{i}),
\end{equation}
the sequence $(X_{n},\xi_{n})_{n\in\mathbf Z_{+}}$ is properly weighted with respect to $\pi$. It is immediately seen that the associated jump process is the original MH output $(Y_{t})_{t\in\mathbf Z_{+}}$ which is known to be a pure Markov chain (rather than a general discrete time semi--Markov process). 

The above analysis suggests that we are allowed to use any distribution (beyond the geometric) for the weights provided (\ref{MH E xi}) is satisfied. However, direct calculation of the importance weight is in general computationally costly or even infeasible making such a task hard to accomplish. Moreover, the geometric distribution comes out naturally, since each simulation from $g(\cdot|\cdot)$ automatically generates the corresponding geometric weight. 

An important issue is the rate of convergence of the resulting chain to the target distribution. In the case of independence MH, i.e.\ when $q(y|z)=q(y)$, Mengersen and Tweedie (1996) showed that if $\tilde{w}^{*}=\sup_{x}\pi(x)/q(x)<\infty$, the algorithm is uniformly ergodic with $\|\mathbf P\{Y_{t}\in\cdot\}-\pi\|\leqslant (1-1/\tilde{w}^{*})^{t}$, $t\in\mathbf Z_{+}$, where $\|\mu\|=\sup_{A\in\mathcal B(\mathcal X)}|\mu(A)|$ denotes as usual the total variation of the signed measure $\mu$. In general, this is the case if $\sup_{x,z}\pi(x)/q(x|z)<\infty$. 

\subsection{Independence samplers with geometric weights}

With the term independence sampler we refer to a simulation scheme where the generated (unweighted) random variates are mutually independent. For instance, this excludes the independence MH algorithm. 

Sahu and Zhigljavsky (2003) and G\r{a}semyr (2002) proposed independence samplers with geometric weighting distributions resulting in ordinary Markov chains. The independence sampler of Sahu and Zhigljavsky (2003) can be described as follows. Let $Z=(Z_{n})_{n\in\mathbf N}$ be a sequence consisting of iid random variates from some distribution $\tilde{g}$. To each $Z_{n}$ it is associated a weight $\xi_{n}$ drawn from the geometric distribution with probability mass function
\[
\mathbf P(\xi_{n}=m|Z_{n}=z) = \frac{1}{1+\kappa \tilde{w}(z)} \left\{ \frac{\kappa \tilde{w}(z)}{1+\kappa \tilde{w}(z)}\right\}^{m},\ m=0,1,2,\ldots,
\]
where $\tilde{w}(x):=\pi(x)/\tilde{g}(x)$. When $\xi_{n}=0$ then the corresponding $Z_{n}$ is rejected. Let $X=(X_{n})_{n\in\mathbf N}$ be the sequence of the {\em accepted} $Z_{n}$'s, i.e.\ those having been weighted by $\xi_{n}>0$. Clearly, the sequence $X$ also consists of iid draws but from the distribution
\[
g(x) = \frac{\tilde{g}(x)P(\xi>0|x)}{\int_{\mathcal X}\tilde{g}(z)P(\xi>0|z)\mu({\rm d}z)} = 
\frac{\pi(x)/[1+\kappa\tilde{w}(x)]}{\int_{\mathcal X}\pi(z)/[1+\kappa\tilde{w}(z)]\mu({\rm d}z)}. 
\]
Moreover, $X_{n}$ is weighted by $\xi_{n}$ which is generated from the (truncated) geometric distribution with probability mass function
\[
\mathbf P(\xi_{n}=m|X_{n}=x) = \frac{1}{1+\kappa \tilde{w}(x)} \left\{ \frac{\kappa \tilde{w}(x)}{1+\kappa \tilde{w}(x)}\right\}^{m-1},\ m=1,2,\ldots.
\]
Since 
\[
\mathbf E\{\xi_{n}|X_{n}=x\} = 1+\kappa\tilde{w}(x) \propto 
w(x) := \frac{\pi(x)}{g(x)} = \left\{\int_{\mathcal X}\frac{\pi(z)}{1+\kappa\tilde{w}(z)}\:\mu({\rm d}z)\right\} \{1+\kappa\tilde{w}(x)\}
\]
the sequence $(X_{n},\xi_{n})_{n\in\mathbf N}$ is properly weighted with respect to $\pi$ and thus the associated (discrete time Markov) jump process $Y=(Y_{t})_{t\in\mathbf Z_{+}}$ converges to $\pi$. Moreover, if $\tilde{w}^{*}=\sup_{x}\pi(x)/\tilde{g}(x)<\infty$, the total variation distance between $\mathbf P\{Y_{t}\in\cdot\}$ and $\pi$ is no greater than $(1+\kappa\tilde{w}^{*})^{-t}$. 

G\r{a}semyr (2002) generalised the above sampler by modifying the rejection rule. After $Z_{n}\sim\tilde{g}$ is drawn, G\r{a}semyr accepts or rejects it according to the result of a Bernoulli trial with some probability of success $q(z_{n})$. Provided it is accepted, it is weighted by a geometric random variate with probability mass function
\[
\mathbf P(\xi_{n}=m|Z_{n}=z) = a(z)\{1-a(z)\}^{m-1},\ m=1,2,\ldots,
\]
where $a(z)\propto q(z)/\tilde{w}(z)$. Taking in particular $a(z)=\{1+\kappa\tilde{w}(z)\}^{-1}$ and $q(z) = \kappa\tilde{w}(z)\{1+\kappa\tilde{w}(z)\}^{-1}$ the sampler reduces to that of Sahu and Zhigljavsky (2003). G\r{a}semyr (2002) shows that the choice $q(z)=\min\{1,\kappa\tilde{w}(z)\}$ and $a(z)=\min\{1,1/\kappa\tilde{w}(z)\}$ is optimal in the sense that it minimizes the asymptotic variance of the estimators $\hat{h}_{n}$ in (\ref{xi estimate}). 

As before, let $X=(X_{n})_{n\in\mathbf N}$ be the sequence of the accepted $Z$'s. For the above optimal choice of $q$, the $X_{n}$'s are independent draws from the distribution 
\[
g(x) = \frac{\tilde{g}(x)\min\{1,\kappa\tilde{w}(x)\}}{\int_{\mathcal X}\tilde{g}(z)\min\{1,\kappa\tilde{w}(z)\}\mu({\rm d}z)} = 
\frac{\pi(x)\min\{1,1/\kappa\tilde{w}(x)\}}{\int_{\mathcal X}\pi(z)\min\{1,1/\kappa\tilde{w}(z)\}\mu({\rm d}z)} 
\]
and
\[
\mathbf E\{\xi_{n}|X_{n}=x\} 
\propto w(x) := \frac{\pi(x)}{g(x)} = \left(
{\int_{\mathcal X}\pi(z)\min\left\{1,\textstyle\frac{1}{\kappa\tilde{w}(z)}\right\}\mu({\rm d}z)} \right) \max\{1,\kappa\tilde{w}(x)\}.
\]
Hence, the sequence $(X_{n},\xi_{n})_{n\in\mathbf N}$ is properly weighted with respect to $\pi$ and thus, the associated (discrete time Markov) jump process $Y=(Y_{t})_{t\in\mathbf Z_{+}}$ converges to $\pi$. Once again, when $\tilde{w}^{*}<\infty$, this Markov chain is uniformly ergodic. When $\kappa>1/\tilde{w}^{*}$, an upper bound of the total variation distance between $\mathbf P\{Y_{t}\in\cdot\}$ and $\pi$ is $(1-1/\kappa\tilde{w}^{*})^{t}$. The case $\kappa\leqslant 1/\tilde{w}^{*}$ correspondes to rejection sampling and consequently $Y_{t}\sim\pi$, $\forall\ t\in\mathbf Z_{+}$.

\subsection{Exponential weights: Pure Markov jump processes}

Consider the case where conditional on $X_{n}=x$, $\xi_{n}$ follows an exponential distribution with mean $\kappa w(x)$. In this case, the associated jump process is a continuous time pure Markov jump process. A particular case of such a sampling scheme is the birth and death MCMC algorithm (Stephens, 2000, Capp\'{e} et al., 2003). 

When $w(x)$ is bounded above, the process exhibits a uniformly ergodic behavior. We state this assertion when $X$ is a Doeblin chain in the following proposition the proof of which can be found in the Appendix.

\begin{proposition} \label{exponential bound}
Let $(X_{n},\xi_{n})_{n\in\mathbf Z_{+}}$ be a weighted sequence satisfying the following:
\begin{enumerate}\setlength{\itemsep}{-.5ex}\vspace*{-1ex} 
\item The sequence $X=(X_{n})_{n\in\mathbf Z_{+}}$ is an ergodic Markov chain with stationary distribution $g$ and its transition distribution $g(\cdot|\cdot)$ satisfies a Doeblin condition, equivalently, there exists a probability distribution $g_{0}$ and a nonnegative constant $\beta$ such that
\begin{equation} \label{Doeblin 1}
g(z|y) \geqslant \beta g_{0}(z),\ \forall\:z,y\in\mathcal X.
\end{equation}
\item The conditional distribution of $\xi_{n}$ given $X_{n}=x$ is exponential with mean $\kappa w(x)=\kappa\pi(x)/g(x)$, where $\pi$ is a probability distribution, and $w^{*}=\sup_{x\in\mathcal X} w(x) <\infty$. 
\end{enumerate}
Then, for the associated Markov jump process $Y=(Y_{t})_{t\geqslant 0}$ it holds
\begin{equation} \label{total variation bound}
\|\mathbf P(Y_{t}\in \cdot)-\pi\| \leqslant \exp\{-\beta t/\kappa w^{*}\},\ t\geqslant 0.
\end{equation}
\end{proposition}

In case the $X_{n}$'s are independent $g$--distributed random variates, (\ref{Doeblin 1}) holds with $\beta=1$ and $g_{0}=g$, thus the total variation distance in (\ref{total variation bound}) is no greater than $\exp\{-t/\kappa w^{*}\}$. This bound is comparable to that of the previous cases, although the actual importance distribution is different (since here it has not been taken care for rejection control).  

\begin{example}  \label{exponential normal mixture}
{\rm
Consider again the distributions of Example \ref{IS normal mixture} for which it can be verified that $w^{*}\approx 6.905$. Let $X=(X_{n})_{n\in\mathbf Z_{+}}$ be iid Cauchy $\mathcal C(0,10)$ random variates and $\xi_{n}$ be exponential with mean $w(x_{n})$. Then, by Proposition \ref{exponential bound}, for the associated Markov jump process $Y$ it holds $\|\mathbf P(Y_{t}\in \cdot)-\pi\| \leqslant \exp\{-t/6.905\}$, $t\geqslant 0$. For instance, for $t\geqslant 31.8$, the total variation distance will be less than $0.01$. 
}
\end{example}

\subsection{Importance weighting an MCMC output}

Let $y_{0},y_{1},y_{2},\ldots$ be the output of any MCMC updating scheme having target distribution $\pi$ and updating distribution $\pi(y_{n}|y_{n-1})$. A crude method for diagnosing convergence of the chain to $\pi$ is checking the convergence of many estimators of some quantities of interest $\mathbf E_{\pi}(h)$ (cf.\ Robert and Casella, 1999, p.382). 

One particular estimator used in this context is the IS estimator
\[
\tilde{h}_{n} := \left.\sum_{i=1}^{n}\frac{\pi(y_{i})}{\pi(y_{i}|y_{i-1})}\:h(y_{i})\right/\sum_{i=1}^{n}\frac{\pi(y_{i})}{\pi(y_{i}|y_{i-1})}
\]
provided that the ratios $\pi(y_{i})/\pi(y_{i}|y_{i-1})$ can be calculated up to a constant. By setting $\xi_{i}:=\pi(y_{i})/\pi(y_{i}|y_{i-1})$, it can be seen that the jump process associated with the weighted sequence $(y_{i},\xi_{i})_{i\in\mathbf N}$ also converges to $\pi$. Indeed, let $x_{i}=(x_{i}^{(1)},x_{i}^{(2)}):=(y_{i-1},y_{i})$. Then, the sequence $(x_{n})_{n\in\mathbf N}$ forms a Markov chain with transition density 
\[
g(x_{i}|x_{i-1}) = 
\delta_{x_{i-1}^{(2)}}(x_{i}^{(1)})
\pi(x_{i}^{(2)}|x_{i}^{(1)}),
\]
where $\delta_{x}(\cdot)$ denotes the Dirac measure at $x$, and limit distribution $g(x)=\pi(x^{(1)})\pi(x^{(2)}|x^{(1)})$. Setting
\[
\xi_{i} = \frac{\pi(x_{i}^{(1)})\pi(x_{i}^{(2)})}{\pi(x_{i}^{(1)})\pi(x_{i}^{(2)}|x_{i}^{(1)})} = 
\frac{\pi(y_{i})}{\pi(y_{i}|y_{i-1})},\ i=1,2,\ldots,
\]
we are led to the sequence $(x_{i},\xi_{i})_{i\in\mathbf N}$ which is properly weighted with respect to the product $\pi(x^{(1)})\pi(x^{(2)})$. Thus, the jump process associated with the (marginal) sequence $(x^{(2)}_{i},\xi_{i})_{i\in\mathbf N}=(y_{i},\xi_{i})_{i\in\mathbf N}$ has $\pi$ as limit distribution. 

Robert and Casella (1999) discuss another approach providing convergent estimators in the context of MH algorithms. As in Subsection 3.1, consider a general MH algorithm with target distribution $\pi$ and proposal $q(\cdot|\cdot)$. Denote now by $y_{0},y_{1},\ldots,y_{M}$ the {\em whole} simulated output (containing also the rejected samples) and by $y_{0},y_{\sigma_{1}},\ldots,y_{\sigma_{n}}$ the {\em accepted} variates. Define $\tau_{i}:=\sup\{\sigma_{j}:\:\sigma_{j}<i\}$ and notice that $y_{\tau_{i}}$ is the last accepted value before time $i$ implying that $y_{i}$ has been generated from $q(\cdot|y_{\tau_{i}})$. The IS estimator of $\mathbf E_{\pi}(h)$ used in this context is 
\[
\tilde{h}_{n} := \left.\sum_{i=1}^{n}\frac{\pi(y_{i})}{q(y_{i}|y_{\tau_{i}})}\:h(y_{i})\right/\sum_{i=1}^{n}\frac{\pi(y_{i})}{q(y_{i}|y_{\tau_{i}})}.
\]
Setting $\xi_{i}=\pi(y_{i})/q(y_{i}|y_{\tau_{i}})$, it can be seen that the jump process associated with the weighted sequence $(y_{i},\xi_{i})_{i\in\mathbf N}$ converges to $\pi$. To see this, define $x_{i}=(x_{i}^{(1)},x_{i}^{(2)}):=(y_{\tau_{i}},y_{i})$. If $\tau_{i}<i-1$, then $\tau_{i}=\tau_{i-1}$ and thus $x_{i}^{(1)}=x_{i-1}^{(1)}$. On the other hand, $\tau_{i}=i-1$ implies that $x_{i}^{(1)}=x_{i-1}^{(2)}$. The latter occurs with probability $a(x_{i-1}^{(1)},x_{i-1}^{(2)})$ (given in (\ref{MH acceptance probability})), 
whereas the former with probability $1-a(x_{i-1}^{(1)},x_{i-1}^{(2)})$. Hence the sequence $x_{1},x_{2},\ldots$ forms a Markov chain with transition density
\[
g(x_{i}|x_{i-1}) = \left\{[1-a(x_{i-1}^{(1)},x_{i-1}^{(2)})]\:\delta_{x_{i-1}^{(1)}}(x_{i}^{(1)})
+ a(x_{i-1}^{(1)},x_{i-1}^{(2)})\:\delta_{x_{i-1}^{(2)}}(x_{i}^{(1)})
\right\} q(x_{i}^{(2)}|x_{i}^{(1)}).
\]
The stationary distribution of this chain is $\pi(x^{(1)})q(x^{(2)}|x^{(1)})$ since $y_{0},y_{\tau_{1}},y_{\tau_{2}},\ldots$ is the original MH output. Weighting each $x_{i}$ by
\[
\xi_{i} = \frac{\pi(x_{i}^{(2)})\pi(x_{i}^{(1)})}{\pi(x_{i}^{(1)})q(x_{i}^{(2)}|x_{i}^{(1)})} = \frac{\pi(y_{i})}{q(y_{i}|y_{\tau_{i}})}
\]
we obtain the desired result for the marginal sequence $(x_{i}^{(2)},\xi_{i})_{i\in\mathbf N}=(y_{i},\xi_{i})_{i\in\mathbf N}$. 


\section{Stationary weighted samples}
\label{section: stationarity samples}

Hereafter we denote by $p(v|x)\nu({\rm d}v)$ the conditional distribution of $\xi$ given $x$, with $\nu$ denoting either the Lebesgue measure (continuous weights) or the counting measure (discrete weights). Moreover, we set 
\[
\overline{P}(u|x) := \int_{[u,\infty)}p(v|x)\nu({\rm d}v) = \mathbf P(\xi\geqslant u|x). 
\]
In this section, we will discuss a slight modification of the standard weighting method under which one can obtain stationary weighted samples from the target distribution $\pi$. We will need first some facts from the theory of semi--Markov processes.

Let $(Y_{t})_{t\geqslant 0}$ be a semi--Markov process with embedded Markov chain $(X_{n})_{n\in\mathbf Z_{+}}$ and respective sojourn times $(\xi_{n})_{n\in\mathbf Z_{+}}$. Let also $S_{n}$ and $N_{t}$ be as in Definition \ref{properly weighted sequence}. Then, $(X_{n},S_{n})_{n\in\mathbf Z_{+}}$ is a Markov renewal process and $(N_{t})_{t\geqslant 0}$ is its counting process. The corresponding excess life (residual age) process is defined by $(V_{t})_{t\geqslant 0}$ by $V_{t}:=S_{N_{t}+1}-t$. Then, $(Y_{t},V_{t})_{t\geqslant 0}$ is a Markov process with stationary distribution $\pi(y)p_{\rm e}(v|y)$, where 
\[
p_{\rm e}(v|y) := \frac{\overline{P}(v|y)}{\kappa w(y)} 
\]
(cf.\ McDonald, 1977).

In the context of weighted samples this suggests the following. Assume that it is possible to generate $X_{0}\sim\pi$. Then, if $X_{0}$ is weighted by $\xi_{0}$ drawn from $p_{\rm e}(\cdot|x_{0})$ (instead of $p(\cdot|x_{0})$), the excess life process $(Y_{t},V_{t})_{t\geqslant 0}$ starts in equilibrium and thus $Y_{t}\sim\pi$, $\forall\:t\geqslant 0$. As a consequence, the estimator
\[
\hat{h}_{t} := t^{-1}\int_{0}^{t}h(Y_{s})\nu({\rm d}s) = t^{-1}\textstyle\left\{\sum_{j=0}^{N_{t}-1}\xi_{j}h(X_{j})+(t-S_{N_{t}})h(X_{N_{t}})\right\}
\]
is exactly unbiased for $\mathbf E_{\pi}(h)$ for any {\em fixed} time $t>0$.

In the case of IS we have that $p_{\rm e}(v|y)\sim \mathcal U(0,w(y))$. Hence, if $X_{0}\sim\pi$, $U\sim\mathcal U(0,1)$ and $X_{1},X_{2},\ldots$ are iid draws from $g$ (or an ergodic Markov chain with stationary distribution $g$), the estimator $\hat{h}_{t}$ becomes 
\[
\hat{h}_{t} = t^{-1}\textstyle\left\{Uw(X_{0})h(X_{0})+\sum_{j=1}^{N_{t}-1}w(X_{j})h(X_{j})+(t-S_{N_{t}})h(X_{N_{t}})\right\}.
\]

In general, simulation from $p_{\rm e}(v|y)$ can be completed in two stages. Generate first a random variate $T\sim q(t|y) \propto t p(t|y)$ and then draw $V$ uniformly distributed in $(0,T)$ (continuous case) or in $\{1,\ldots,T\}$ (discrete case). This is a quite easy task for many standard distributions. For example, if $p$ is a gamma distribution with shape $a$, then $q$ is also a gamma distribution but with shape $a+1$. Or, if $p$ is Poisson then $q$ is a truncated Poisson with the same parameter. 

When $(X_{n})_{n\in\mathbf Z_{+}}$ is an iid sequence from $g$ and the hazard rate $p(v|y)/\overline{P}(v|y)$ is uniformly bounded away from zero, i.e.\ if
\begin{equation} \label{bounded hazard}
\varepsilon_{*} := \inf_{v,y}\frac{p(v|y)}{\overline{P}(v|y)} > 0,
\end{equation}
an accept--reject argument can give at once a stationary semi--Markov process. Indeed, let $(x,\xi)$ be a draw from $g(x)p(\xi|x)$. Then, 
\[
\frac{\pi(x)p_{e}(\xi|x)}{g(x)p(\xi|x)} = \frac{\overline{P}(\xi|x)}{\kappa p(\xi|x)} \leqslant \frac{1}{\kappa\varepsilon_{*}} 
\]
and consequently if $U$ is an independent $\mathcal U(0,1)$ random variate and $U\leqslant \varepsilon_{*}\overline{P}(\xi|x)/p(\xi|x)$ holds, $(x,\xi)$ can be considered as $\pi(x)p_{e}(\xi|x)$--distributed. 

\begin{lemma} \label{lemma:bounded hazard}
If (\ref{bounded hazard}) holds then (i) $w^{*}=\sup_{x\in\mathcal X}w(x)<\infty$ and (ii) there exists $a>1$ such that $\mathbf E_{g}\{a^{\xi}\}:=\int\int a^{v}g(y)p(v|y)\nu({\rm d}v)\mu({\rm d}y)<\infty$. 
\end{lemma}

According to the part (i) of this lemma, the above approach does not offer much, since the accept--reject method applies also to the $x_{n}$'s. However, part (ii) can be used to obtain a crude bound for total variation distance to stationarity for the associated semi--Markov process.

\begin{theorem} \label{Yt sim pi}
Assume that $(X_{n})_{n\in\mathbf Z_{+}}$ is an iid sequence from $g$ and (\ref{bounded hazard}) holds. Then for the associated jump process $Y=(Y_{t})_{t\geqslant 0}$ we have the following: \\
\hspace*{2ex}{\rm (a)} There is an almost surely finite time $\tau\geqslant 0$, such that $Y_{t}\sim\pi$ for $t\geqslant\tau$. \\
\hspace*{2ex}{\rm (b)} The total variation distance between the law of $Y_{t}$ and $\pi$ converges to zero exponentially fast in $t$.
\end{theorem}


\section*{Appendix}
\appendix
\setcounter{section}{1}

In order to prove Proposition \ref{exponential bound} we will need the following lemma. 

\noindent
\begin{lemma} Under the assumptions of Proposition \ref{exponential bound} it holds
\begin{equation} \label{eq:exponential limit}
\lim_{t\downarrow 0} \left\{1-\beta\big(1-e^{-t/\kappa w^{*}}\big)\int_{\mathcal X}e^{-t/\kappa w(z)} g_{0}(z)\mu({\rm d}z)\right\}^{1/t}  = \exp\{-\beta/\kappa w^{*}\}. 
\end{equation}
\end{lemma}

\begin{proof} Write first
\begin{eqnarray*}
a(t) &:=& 1-\beta\big(1-e^{-t/\kappa w^{*}}\big)\int_{\mathcal X}e^{-t/\kappa w(z)} g_{0}(z)\mu({\rm d}z)
\\
&=&
1-\beta\big(1-e^{-t/\kappa w^{*}}\big) + \beta\big(1-e^{-t/\kappa w^{*}}\big)\int_{\mathcal X}\big(1-e^{-t/\kappa w(z)}\big) g_{0}(z)\mu({\rm d}z)  
\end{eqnarray*}
Observe that (\ref{Doeblin 1}) implies $g(z)\geqslant \beta g_{0}(z)$, $\forall\:z$. Hence, 
\[
\beta\int_{\mathcal X}\big(1-e^{-t/\kappa w(z)}\big) g_{0}(z)\mu({\rm d}z) \leqslant \int_{\mathcal X}\big(1-e^{-t/\kappa w(z)}\big) g(z)\mu({\rm d}z) = \mathbf P_{g}(\xi_{0}\leqslant t) = o(1).  
\]
Since $1-e^{-t/\kappa w^{*}}=O(t)$, we have that $\beta\big(1-e^{-t/\kappa w^{*}}\big)\int_{\mathcal X}\big(1-e^{-t/\kappa w(z)}\big) g_{0}(z)\mu({\rm d}z) = o(t)$ and thus,
\[ 
1-\beta\big(1-e^{-t/\kappa w^{*}}\big)\leqslant a(t) \leqslant 1-\beta\big(1-e^{-t/\kappa w^{*}}\big)+o(t),
\]
Raising to the power of $1/t$ and taking the limits we obtain the desired result. 
\end{proof}

\vspace*{1ex}

\noindent
{\it Proof of Proposition \ref{exponential bound}.} 
Since $(Y_{t})_{t\geqslant 0}$ is a Markov process, every $\delta$--skeleton, i.e.\ any sequence $(Y_{n\delta})_{n\in\mathbf Z_{+}}$ for fixed $\delta>0$, forms a Markov chain with transition kernel $\mathbf P^{\delta}(y,A)=\mathbf P\{Y_{\delta}\in A|Y_{0}=y\}$. But then,
\begin{eqnarray*}
\mathbf P^{\delta}(y,A)
&=& 
\mathbf P\{\xi_{0}>\delta,y\in A|X_{0}=y\} + \sum_{m=1}^{\infty}\mathbf P\left\{\left.\textstyle S_{m}\leqslant \delta <S_{m+1},X_{m}\in A\right|X_{0}=y\right\}
\\
&\geqslant&
\mathbf P\{\xi_{0}\leqslant \delta,\xi_{1}>\delta,X_{1}\in A|X_{0}=y\}
\\
&=&
\big\{1-e^{-\delta/\kappa w(y)}\big\} \int_{A}e^{-\delta/\kappa w(z)} g(z|y)\mu({\rm d}z)
\\
&\geqslant&
\big\{1-e^{-\delta/\kappa w^{*}}\big\} \int_{A}e^{-\delta/\kappa w(z)} \beta g_{0}(z)\mu({\rm d}z) 
\\
&=&
\varepsilon_{\delta}Q_{\delta}(A),\ \forall\:y,
\end{eqnarray*}
where 
\[
\varepsilon_{\delta}:=\beta\big\{1-e^{-\delta/\kappa w^{*}}\big\} 
\int_{\mathcal X}e^{-\delta/\kappa w(z)} g_{0}(z)\mu({\rm d}z)
\]
and 
\[
Q_{\delta}(A) := \frac{\int_{A}e^{-\delta/\kappa w(z)} g_{0}(z)\mu({\rm d}z)}{\int_{\mathcal X}e^{-\delta/\kappa w(z)} g_{0}(z)\mu({\rm d}z)}.
\]
Thus, the $\delta$--skeleton chain satisfies a Doeblin condition. It is well known that this implies  
\[
\|\mathbf P(Y_{n\delta}\in\cdot)-\pi\| \leqslant (1-\varepsilon_{\delta})^{n},\ \forall\:n\in\mathbf Z_{+},
\]
hence, writing $t=n\delta$, we conclude that
\[
\|\mathbf P(Y_{t}\in\cdot)-\pi\| \leqslant \{(1-\varepsilon_{t/n})^{n/t}\}^{t},\ \forall\:n\in\mathbf Z_{+}. 
\]
However, $\displaystyle\lim_{n\uparrow\infty}(1-\varepsilon_{t/n})^{n/t}=\displaystyle\lim_{t\downarrow 0}(1-\varepsilon_{t})^{1/t}$ and the result follows from Lemma A.1.

\vspace*{3ex}

\noindent
{\it Proof of Lemma \ref{lemma:bounded hazard}.}  
(i) We have that
\begin{eqnarray*}
&&\frac{p(v|y)}{\bar{P}(v|y)}>\varepsilon_{*},\ \forall\:v,y
\ \Leftrightarrow\ p(v|y)>\varepsilon_{*}\overline{P}(v|y),\ \forall\:v,y
\\ 
&\Rightarrow& \int p(v|y)\nu({\rm d}v) > \varepsilon_{*}\int \overline{P}(v|y) \nu({\rm d}v),\ \forall\:y
\ \Leftrightarrow\ 1>\varepsilon_{*}\kappa w(y),\ \forall\:y
\\
&\Leftrightarrow& w(y)<1/\kappa\varepsilon_{*},\ \forall\:y,
\end{eqnarray*}
thus $w^{*}=\sup_{y}w(y)\leqslant 1/\kappa\varepsilon_{*}<\infty$. 

\noindent
(ii) We will first prove that $\mathbf E\{\xi^{m}|y\}\leqslant m!/\varepsilon_{*}^{m}$, $m\geqslant 1$, by induction. By (i) the assertion holds for $m=1$. Assuming that it holds for some $m$ we have
\begin{eqnarray*}
\mathbf E\{\xi^{m}|y\} &=& \int v^{m}p(v|y)\nu({\rm d}v)\\ &\geqslant& \varepsilon_{*}\int v^{m}\overline{P}(v|y)\nu({\rm d}v) \ = \  
\varepsilon_{*}\int v^{m}\left(\int_{[v,\infty)}p(u|y)\nu({\rm d}u)\right)\nu({\rm d}v)\\ &=& \varepsilon_{*}\int\left(\int_{(0,u]}v^{m}\nu({\rm d}v)\right)p(u|y)\nu({\rm d}u)
\ \geqslant\ \varepsilon_{*}\int\left(\int_{v=0}^{u}v^{m}{\rm d}v\right)p(u|y)\nu({\rm d}u)\\ &=& \varepsilon_{*}\int\frac{u^{m+1}}{m+1}p(u|y)\nu({\rm d}u)\ =\   \frac{\varepsilon_{*}}{m+1}\:\mathbf E\{\xi^{m+1}|y\},\ \forall\:y,
\end{eqnarray*}
and thus, $\mathbf E\{\xi^{m+1}|y\}\leqslant (m+1)!/\varepsilon_{*}^{m+1}$. Now,
\[
\mathbf E\{a^{\xi}|y\}= \mathbf E\{e^{\xi\log a}|y\}=\sum_{m=0}^{\infty}\frac{(\log a)^{m}}{m!}\:\mathbf E\{\xi^{m}|y\} \leqslant \sum_{m=0}^{\infty}\frac{(\log a)^{m}}{m!}\:\frac{m!}{\varepsilon_{*}^{m}} = \frac{1}{1-\log a/\varepsilon_{*}}<\infty
\]
for any $a\in[1,\exp\{\varepsilon_{*}\})$, the same holding true also for $\mathbf E_{g}\{a^{\xi}\}=\mathbf E_{g}[\mathbf E\{a^{\xi}|Y\}]$.

\vspace*{3ex}

\noindent
{\it Proof of Proposition \ref{Yt sim pi}.} After $(X_{n},\xi_{n})$ has been drawn, generate an independent random variate $U_{n}\sim\mathcal U(0,1)$. Let $N$ be the first index $n\in\mathbf Z_{+}$ at which $U_{n}\leqslant \varepsilon_{*}\overline{P}(\xi_{n}|x_{n})/p(\xi_{n}|x_{n})$ occurs. Then $(Y_{S_{N}},V_{S_{N}})=(X_{N},\xi_{N})\sim\pi(\cdot)p_{e}(\cdot|\cdot)$. Since $\pi p_{e}$ is stationary for $(Y_{t},V_{t})_{t\geqslant 0}$, we have that $Y_{t}\sim\pi$ for every $t\geqslant \tau:=S_{N}$. Now, by the standard theory of rejection sampling, $N+1\sim\mathcal G(\kappa\varepsilon_{*})$. Moreover, the ``rejected'' $(X,\xi)$'s have the residual distribution 
\[
r(x,v):=
\frac{g(x)p(v|x)-\kappa\varepsilon_{*}\pi(x)p_{e}(v|x)}{1-\kappa\varepsilon_{*}} = 
g(x)\:\frac{p(v|x)-\varepsilon_{*}\overline{P}(v|x)}{1-\kappa\varepsilon_{*}}. 
\]
Notice that for any $a$ such that $\varphi_{g}(a):=\mathbf E_{g}\{a^{\xi}\}<\infty$ it also holds $\varphi_{r}(a):=\mathbf E_{r}\{a^{\xi}\}<\infty$ since $r(x,v)\leqslant g(x)p(v|x)/(1-\kappa\varepsilon_{*})$ for all $x, v$. The function $\varphi_{r}(a)$ is continuous and increasing in $a\geqslant 1$ with $\varphi_{r}(1)=1$, hence there exists $\rho>1$ such that $\varphi_{r}(\rho)(1-\kappa\varepsilon_{*})<1$. Then,
\begin{equation} \label{tau inequality}
\mathbf P(\tau>t) = \mathbf P(S_{N}>t) \leqslant \rho^{-t}\mathbf E\{\rho^{S_{N}}\} = \rho^{-t}\mathbf E\{\varphi_{r}(\rho)^{N}\} =  \frac{\kappa\varepsilon_{*}\rho^{-t}}{1-\varphi_{r}(\rho)+\kappa\varepsilon_{*}\varphi_{r}(\rho)} 
\end{equation}
which converges to zero as $t\rightarrow\infty$ implying $\mathbf P(\tau<\infty)=1$. 

\noindent
(b) Let $A\in\mathcal B(\mathcal X)$. Then, 
\begin{eqnarray*}
|\mathbf P\{Y_{t}\in A\}-\pi(A)| 
&=& 
|\mathbf P\{Y_{t}\in A|\tau>t\}\mathbf P\{\tau>t\} +\mathbf P\{Y_{t}\in A|\tau\leqslant t\}\mathbf P\{\tau\leqslant t\}-\pi(A)| 
\\
&\leqslant&
\mathbf P\{\tau>t\}\:|\mathbf P\{Y_{t}\in A|\tau>t\}-\pi(A)| 
\\
&\leqslant&
\mathbf P\{\tau>t\},
\end{eqnarray*}
since $\mathbf P\{Y_{t}\in A|\tau\leqslant t\}=\pi(A)$ by (a), and $|\mathbf P\{Y_{t}\in A|\tau>t\}-\pi(A)|\leqslant 1$. Thus, 
\[
\|\mathbf P\{Y_{t}\in \cdot\} - \pi\| = \sup_{A\in\mathcal B(\mathcal X)}|\mathbf P\{Y_{t}\in A\}-\pi(A)|\leqslant \mathbf P\{\tau>t\} = O(\rho^{-t})
\] 
by \eqref{tau inequality}.


\section*{References}

\begin{small}
\begin{description}
\item Capp\'e O., Robert C.P. and Ryd\'en T. (2003). Reversible jump, birth--and--death and more general continuous time Markov chain Monte Carlo samplers. Journal of the Royal Statistical Society, Series B 65: 679--700. 
\item G\r{a}semyr J. 2002. Markov chain Monte Carlo algorithms with independent proposal distribution and their relationship to importance sampling and rejection sampling. Technical report, Department of Mathematics, University of Oslo.
\item Hastings W.K. 1970. Monte Carlo sampling methods using Markov chains and their applications. Biometrika 57: 97--109. 
\item Liu J.S. 2001. Monte Carlo Strategies in Scientific Computing. Springer--Verlag, New York.
\item Liu J.S. and Chen R. 1998. Sequential Monte Carlo methods for dynamic systems. Journal of the American Statistical Association 93: 1032--1044.
\item Limnios N. and Opri\c{s}an G. 2001. Semi--Markov processes and reliability. Statistics for industry and technology, Birkh\"{a}user. 
\item Marshall A.W. 1956. The use of multi--stage sampling schemes in Monte Carlo computations. In: Meyer M.A. (Ed.), Symposium on Monte Carlo methods. Wiley, New York, pp. 123--140. 
\item McDonald D.R. 1977. Equilibrium measures for semi-Markov processes. The Annals of Probability 5: 818--822.
\item Metropolis N., Rosenbluth A.W., Rosenbluth M.N., Teller A.H. and Teller E. 1953. Equations of state calculations by fast computing machines. The Journal of Chemical Physics 21: 1087--1091. 
\item Robert C.P. and Casella G. 1999. Monte Carlo statistical methods. Springer--Verlag, New York. 
\item Rubin D.B. 1987. A noniterative sampling/importance resampling alternative to the data augmentation algorithm for creating a few imputations when the fraction of missing information is modest: the SIR algorithm. Journla of the American Statistical Association 82: 543--546.
\item Sahu S.K. and Zhigljavsky A.A. 2003. Self--regenerative Markov chain Monte Carlo with adaptation. Bernoulli 9: 395--422. 
\item Smith A.F.M. and Gelfand A.E. 1992. Bayesian statistics without tears. The American Statistician 46: 84--88.

\end{description}
\end{small}

\end{document}